\documentclass[times,sort&compress,3p,12pt]{elsarticle}
\journal{Results in Mathematics}

\usepackage{amsmath,amsfonts,amssymb,amsthm,graphicx,hyperref,url}
\usepackage[labelfont=bf]{caption}

\usepackage{mathtools}
\mathtoolsset{showonlyrefs=true} 
\usepackage{dsfont} 
\usepackage{mathrsfs} 
\usepackage{subcaption}
\usepackage{xcolor}
\usepackage{geometry}
\theoremstyle{plain}
\newtheorem{theorem}{Theorem}

\newtheorem{corollary}{Corollary}

\theoremstyle{definition}

\newtheorem{remark}{Remark}

\newcommand{\N}{\mathbb{N}}
\newcommand{\R}{\mathbb{R}}

\newcommand{\QQ}{\mathbb{Q}}
\newcommand{\PP}{\mathbb{P}}
\newcommand{\EE}{\mathbb{E}}

\newcommand{\bb}[1]{\boldsymbol{#1}}
\newcommand{\OO}{\mathcal O}

\newcommand{\leqdef}{\vcentcolon=}

\newcommand{\rd}{{\rm d}}
\newcommand{\ind}{\mathds{1}}

\makeatletter \@addtoreset{equation}{section} \makeatother

\allowdisplaybreaks

\begin{document}

\begin{frontmatter}

    \title{On the Le Cam distance between multivariate hypergeometric\\ and multivariate normal experiments}%

    \author[a1,a2]{Fr\'ed\'eric Ouimet}%

    \address[a1]{McGill University, Montreal, QC H3A 0B9, Canada.}%
    \address[a2]{California Institute of Technology, Pasadena, CA 91125, USA.}%

    \ead{frederic.ouimet2@mcgill.ca}%

    \begin{abstract}
        In this short note, we develop a local approximation for the log-ratio of the multivariate hypergeometric probability mass function over the corresponding multinomial probability mass function.
        In conjunction with the bounds from \citet{MR1922539} and \citet{MR4249129} on the total variation between the law of a multinomial vector jittered by a uniform on $(-1/2,1/2)^d$ and the law of the corresponding multivariate normal distribution, the local expansion for the log-ratio is then used to obtain a total variation bound between the law of a multivariate hypergeometric random vector jittered by a uniform on $(-1/2,1/2)^d$ and the law of the corresponding multivariate normal distribution.
        As a corollary, we find an upper bound on the Le Cam distance between multivariate hypergeometric and multivariate normal experiments.
    \end{abstract}

    \begin{keyword}
        multivariate hypergeometric distribution \sep sampling without replacement \sep multinomial distribution \sep normal approximation \sep Gaussian approximation \sep local approximation \sep local limit theorem \sep asymptotic statistics \sep multivariate normal distribution \sep Le Cam distance \sep total variation \sep deficiency \sep comparison of experiments
        \MSC[2020]{Primary: 62E20, 62B15 Secondary: 60F99, 60E05, 62H10, 62H12}
    \end{keyword}

\end{frontmatter}

\section{Introduction}\label{sec:intro}

    Let $d\in \N$.
    The $d$-dimensional (unit) simplex and its interior are defined by
    \begin{equation}\label{eq:def.simplex}
        \mathcal{S}_d \leqdef \big\{\bb{x}\in [0,1]^d: \|\bb{x}\|_1 \leq 1\big\}, \qquad \mathrm{Int}\hspace{0.2mm}(\mathcal{S}_d) \leqdef \big\{\bb{x}\in (0,1)^d: \|\bb{x}\|_1 < 1\big\},
    \end{equation}
    where $\|\bb{x}\|_1 \leqdef \sum_{i=1}^d |x_i|$ denotes the $\ell^1$ norm on $\R^d$.
    Given a set of probability weights $\bb{p}\in N^{-1} \, \N_0^d \cap \mathrm{Int}\hspace{0.2mm}(\mathcal{S}_d)$, the probability mass function of the multivariate hypergeometric distribution,  $\mathrm{Hypergeometric}\hspace{0.2mm}(N,n,\bb{p})$, is defined, by \citet[Chapter~39]{MR1429617}, as
    \vspace{-1mm}
    \begin{equation}\label{eq:hypergeometric.pmf}
        P_{N,n,\bb{p}}(\bb{k}) \leqdef \frac{\prod_{i=1}^{d+1} \binom{N p_i}{k_i}}{\binom{N}{n}}, \quad \bb{k}\in \mathbb{K}_d,
    \end{equation}
    where $p_{d+1} \leqdef 1 - \|\bb{p}\|_1 > 0$, $k_{d+1} \leqdef n - \|\bb{k}\|_1$, $n,N\in \N$ with $n \leq N$, and
    \begin{equation}
        \mathbb{K}_d \leqdef \big\{\bb{k}\in \N_0^d \cap n \mathcal{S}_d : k_i\in [0,N p_i] \, \, \text{for all } 1 \leq i \leq d + 1\big\}.
    \end{equation}
    This distribution represents the first $d$ components of the vector of categorical sample counts when randomly sorting a random sample of $n$ objects from a finite population of $N$ objects into $d+1$ categories, where $p_i, ~1 \leq i \leq d + 1$, is the probability of any given object to be sorted in the $i$-th category.

    \vspace{2mm}
    Our first main goal in this paper is to develop a local approximation for the log-ratio of the multivariate hypergeometric probability mass function \eqref{eq:hypergeometric.pmf} over the $\mathrm{Multinomial}\hspace{0.2mm}(n,\bb{p})$ probability mass function, namely
    \begin{equation}\label{eq:multinomial.pmf}
        Q_{n,\bb{p}}(\bb{k}) \leqdef \frac{n!}{\prod_{i=1}^{d+1} k_i!} \prod_{i=1}^{d+1} p_i^{k_i}, \quad \bb{k}\in \N_0^d \cap n \mathcal{S}_d.
    \end{equation}
    This latter distribution represents the exact same as \eqref{eq:hypergeometric.pmf} above, except that the population from which the $n$ objects are drawn is infinite ($N = \infty$).
    Another way of distinguishing $P_{N,n,\bb{p}}$ and $Q_{n,\bb{p}}$ in a finite population of $N$ objects is to say that we sample the $n$ objects without replacement and with replacement, respectively.
    In both cases, the categorical probabilities $(\bb{p},p_{d+1})$ are the same.
    For good general references on normal approximations, we refer the reader to \citet{MR0436272} and \citet{MR1295242}.

    \vspace{2mm}
    Our second main goal is to prove an upper bound on the total variation between the probability measure on $\R^d$ induced by a random vector distributed according to $P_{N,n,\bb{p}}$ and jittered by a uniform on $(-1/2,1/2)^d$, and the probability measure on $\R^d$ induced by a multivariate normal random vector with the same mean and covariances as a random vector distributed according to $Q_{n,\bb{p}}$, namely $n \bb{p}$ and $n (\mathrm{diag}(\bb{p}) - \bb{p} \bb{p}^{\top})$.
    The proof makes use of the total variation bound from \citet[Lemma~3.1]{MR4249129} (which improved Lemma~2 in \cite{MR1922539}) on the total variation between the probability measure on $\R^d$ induced by a multinomial vector distributed according to $Q_{n,\bb{p}}$ and jittered by a uniform on $(-1/2,1/2)^d$, and the probability measure on $\R^d$ induced by a multivariate normal random vector with the same mean and covariances.
    As pointed out by \citet[p.732]{MR3717995}, the univariate case here would be much simpler since \citet[p.62-63]{MR0254884} showed that the hypergeometric probability mass function can be written as a ratio of three binomial probability mass functions, and local limit theorems are well-known for the binomial distribution, see, e.g., \citet{MR56861} and \citet{MR207011}.

    \vspace{2mm}
    The deficiency between a given statistical experiment and another measures the loss of information from carrying inferences to the second setting using information from the first setting. This loss of information goes in both directions, but the deficiency is not necessarily symmetric. The maximum of the two deficiencies is called the Le Cam distance (or $\Delta$-distance in \cite{MR1784901}).
    The usefulness of this notion comes from the fact that seemingly completely different statistical experiments can result in asymptotically equivalent inferences using Markov kernels to carry information from one setting to another. For instance, it was famously shown by \citet{MR1425959} that the density estimation problem and the Gaussian white noise problem are asymptotically equivalent in the sense that the Le Cam distance between the two experiments goes to $0$ as the number of observations goes to infinity. The main idea was that the information we get from sampling observations from an unknown density function and counting the observations that fall in the various boxes of a fine partition of the density's support can be encoded using the increments of a properly scaled Brownian motion with drift $t\mapsto \int_0^t \hspace{-1.5mm}\sqrt{f(s)}\hspace{0.5mm}\rd s$, and vice versa.
    An alternative (simpler) proof of this asymptotic equivalence was shown by \citet{MR2102503} who combined a Haar wavelet cascade scheme with coupling inequalities relating the binomial and univariate normal distributions at each step (a similar argument was developed previously by \citet{MR1922539} to derive a multinomial/multivariate coupling inequality).
    Not only \citet{MR2102503} streamlined the proof of the asymptotic equivalence originally shown by \citet{MR1425959}, but their results hold for a larger class of densities and the asymptotic equivalence was also extended to Poisson processes.
    Our third main result in the present paper extends the multinomial/multivariate normal comparison from \cite{MR1922539} (revisited and improved by \citet{MR4249129}, who removed the inductive part of the argument) to the multivariate hypergeometric/multivariate normal comparison (recall from \eqref{eq:multinomial.pmf} that the multinomial distribution is just the limiting case $N = \infty$ of the multivariate hypergeometric distribution).
    For an excellent and concise review on Le Cam's theory for the comparison of statistical models, we refer the reader to \citet{MR3850766}.

    \vspace{2mm}
    The three results we have just described are presented in Section~\ref{sec:results}, and the related proofs are gathered in Section~\ref{sec:proofs}.
    Here are now some motivations for these results.
    First, we believe that the first two results (the local expansion of the log-ratio and the total variation bound) could help in developing asymptotic Berry-Esseen type bounds for the symmetric multivariate hypergeometric distribution and the symmetric multinomial distribution, similar to the exact optimal bounds proved recently by \citet{MR3717995} in the univariate setting.
    Second, there might be a way to use the Le Cam distance upper bound between multivariate hypergeometric and multivariate normal experiments to extend the results on the asymptotic equivalence between the density estimation problem and the Gaussian white noise problem shown by \citet{MR1425959} and \citet{MR2102503}.

    \begin{remark}
        Throughout the paper, the notation $u = \OO(v)$ means that $\limsup_{N\to \infty} |u / v| < C$, where $C\in (0,\infty)$ is a universal constant.
        Whenever $C$ might depend on a parameter, we add a subscript (for example, $u = \OO_d(v)$).
    \end{remark}

\section{Results}\label{sec:results}

    Our first main result is an asymptotic expansion for the log-ratio of the multivariate hypergeometric probability mass function \eqref{eq:hypergeometric.pmf} over the corresponding multinomial probability mass function \eqref{eq:multinomial.pmf}.

    \begin{theorem}[Local limit theorem for the log-ratio]\label{thm:p.k.expansion}
        Assume that $n,N\in \N$ with $n\leq N$ and $\bb{p}\in N^{-1} \, \N_0^d \cap \mathrm{Int}\hspace{0.2mm}(\mathcal{S}_d)$ hold, and pick any $\gamma\in (0,1)$.
        Then, uniformly for $\bb{k}\in \mathbb{K}_d$ such that $\max_{1 \leq i \leq d + 1} (k_i / p_i) \leq \gamma N$ and $n \leq \gamma N$, we have, as $N\to \infty$,
        \begin{equation}\label{eq:thm:p.k.expansion.eq.1}
            \log\left(\frac{P_{N,n,\bb{p}}(\bb{k})}{Q_{n,\bb{p}}(\bb{k})}\right) = \frac{1}{N} \left[\bigg(\frac{n^2}{2} - \frac{n}{2}\bigg) - \sum_{i=1}^{d+1} \frac{1}{p_i} \cdot \bigg(\frac{k_i^2}{2} - \frac{k_i}{2}\bigg)\right] + \OO_{\gamma}\left(\frac{1}{N^2} \left[n^3 + \sum_{i=1}^{d+1} \frac{k_i^3}{p_i^2}\right]\right),
        \end{equation}
        and
        \begin{equation}\label{eq:thm:p.k.expansion.eq.2}
            \begin{aligned}
                \log\left(\frac{P_{N,n,\bb{p}}(\bb{k})}{Q_{n,\bb{p}}(\bb{k})}\right)
                &= \frac{1}{N} \left[\bigg(\frac{n^2}{2} - \frac{n}{2}\bigg) - \sum_{i=1}^{d+1} \frac{1}{p_i} \cdot \bigg(\frac{k_i^2}{2} - \frac{k_i}{2}\bigg)\right] \\
                &\quad+ \frac{1}{N^2} \left[\bigg(\frac{n^3}{6} - \frac{n^2}{4} + \frac{n}{12}\bigg) - \sum_{i=1}^{d+1} \frac{1}{p_i^2} \cdot \bigg(\frac{k_i^3}{6} - \frac{k_i^2}{4} + \frac{k_i}{12}\bigg)\right] \\
                &\quad+ \OO_{\gamma}\left(\frac{1}{N^3} \left[n^4 + \sum_{i=1}^{d+1} \frac{k_i^4}{p_i^3}\right]\right).
            \end{aligned}
        \end{equation}
    \end{theorem}

    The local limit theorem above together with the total variation bound in \cite{MR1922539,MR4249129} between jittered multinomials and the corresponding multivariate normals allow us to derive an upper bound on the total variation between the probability measure on $\R^d$ induced by a multivariate hypergeometric random vector jittered by a uniform random vector on $(-1/2,1/2)^d$ and the probability measure on $\R^d$ induced by a multivariate normal random vector with the same mean and covariances as the multinomial distribution in \eqref{eq:multinomial.pmf}.

    \begin{theorem}[Total variation upper bound]\label{thm:total.variation.thm}
        Assume that $n,N\in \N$ with $n\leq (3/4) \, N$ and $\bb{p}\in N^{-1} \, \N_0^d \cap \mathrm{Int}\hspace{0.2mm}(\mathcal{S}_d)$ hold.
        Let $\bb{K}\sim \mathrm{Hypergeometric}\hspace{0.2mm}(N,n,\bb{p})$, $\bb{L}\sim \mathrm{Multinomial}\hspace{0.2mm}(n,\bb{p})$, and $\bb{U}, \bb{V}\sim \mathrm{Uniform}\hspace{0.2mm}(-1/2,1/2)^d$, where $\bb{K}$, $\bb{L}$, $\bb{U}$ and $\bb{V}$ are assumed to be jointly independent.
        Define $\bb{X} \leqdef \bb{K} + \bb{U}$ and $\bb{Y} \leqdef \bb{L} + \bb{V}$, and let $\widetilde{\PP}_{N,n,\bb{p}}$ and $\widetilde{\QQ}_{n,\bb{p}}$ be the laws of $\bb{X}$ and $\bb{Y}$, respectively.
        Also, let $\QQ_{n,\bb{p}}$ be the law of the $\mathrm{Normal}_d(n \bb{p}, n \Sigma_{\bb{p}})$ distribution, where $\Sigma_{\bb{p}} \leqdef \mathrm{diag}(\bb{p}) - \bb{p} \bb{p}^{\top}$.
        Then, as $N\to \infty$,
        \begin{equation}
            \begin{aligned}
                \|\widetilde{\PP}_{N,n,\bb{p}} - \QQ_{n,\bb{p}}\|
                &\leq \sqrt{2 \sum_{i=1}^{d+1} \left(\frac{1}{\nu_i}\right)^{n \nu_i p_i} \left(\frac{1 - p_i}{1 - \nu_i p_i}\right)^{n (1 - \nu_i p_i)} \hspace{-1mm}+ \OO\left(\frac{n^2}{N}\right)} \\
                &\quad+ \OO\left(\frac{d}{\sqrt{n}} \sqrt{\frac{\max\{p_1,\dots,p_d,p_{d+1}\}}{\min\{p_1,\dots,p_d,p_{d+1}\}}}\right),
            \end{aligned}
        \end{equation}
        where $\nu_i \leqdef \lceil p_i^{-1} - 1\rceil$, for $1 \leq i \leq d + 1$, and $\| \cdot \|$ denotes the total variation norm.
    \end{theorem}

    Since the Le Cam distance is a pseudometric and the Markov kernel that jitters a random vector by a uniform on $(-1/2,1/2)^d$ is easily inverted (round off each component of the vector to the nearest integer), then we find, as a consequence of the total variation bound in Theorem~\ref{thm:total.variation.thm}, an upper bound on the Le Cam distance between multivariate hypergeometric and multivariate normal experiments.

    \begin{theorem}[Le Cam distance upper bound]\label{thm:bound.deficiency.distance}
        Assume that $n,N\in \N$ with $n\leq (3/4) \, N$ holds.
        For any given $R\geq 1$, let
        \begin{equation}
            \Theta_R \leqdef \left\{\bb{p}\in N^{-1} \, \N_0^d \cap \mathrm{Int}\hspace{0.2mm}(\mathcal{S}_d) : \, \frac{\max\{p_1,\dots,p_d,p_{d+1}\}}{\min\{p_1,\dots,p_d,p_{d+1}\}} \leq R\right\}.
        \end{equation}
        Define the experiments
        \begin{alignat*}{6}
            &\mathscr{P}
            &&\leqdef &&~\{\PP_{N,n,\bb{p}}\}_{\bb{p}\in \Theta_R}, \quad &&\PP_{N,n,\bb{p}} ~\text{is the measure induced by } \mathrm{Hypergeometric}\hspace{0.2mm}(N,n,\bb{p}), \\
            &\mathscr{Q}\hspace{-0.5mm}
            &&\leqdef &&~\{\QQ_{n,\bb{p}}\}_{\bb{p}\in \Theta_R}, \quad &&\QQ_{n,\bb{p}} ~\text{is the measure induced by } \mathrm{Normal}_d(n \bb{p}, n \Sigma_{\bb{p}}).
        \end{alignat*}
        Then, for $N\geq n^3 / d^{\hspace{0.2mm}2}$, we have the following upper bound on the Le Cam distance $\Delta(\mathscr{P},\mathscr{Q})$ between $\mathscr{P}$ and $\mathscr{Q}$,
        \vspace{-2mm}
        \begin{equation}\label{eq:thm:bound.deficiency.distance.bound}
            \Delta(\mathscr{P},\mathscr{Q}) \leqdef \max\{\delta(\mathscr{P},\mathscr{Q}),\delta(\mathscr{Q},\mathscr{P})\} \leq C_R \, \frac{d}{\sqrt{n}},
        \end{equation}
        where $C_R$ is a positive constant that depends only on $R$,
        \begin{equation}\label{eq:def:deficiency.one.sided}
            \begin{aligned}
                \delta(\mathscr{P},\mathscr{Q})
                &\leqdef \inf_{T_1} \sup_{\bb{p}\in \Theta_R} \bigg\|\int_{\mathbb{K}_d} T_1(\bb{k}, \cdot \, ) \, \PP_{N,n,\bb{p}}(\rd \bb{k}) - \QQ_{n,\bb{p}}\bigg\|, \\
                \delta(\mathscr{Q},\mathscr{P})
                &\leqdef \inf_{T_2} \sup_{\bb{p}\in \Theta_R} \bigg\|\PP_{N,n,\bb{p}} - \int_{\R^d} T_2(\bb{z}, \cdot \, ) \, \QQ_{n,\bb{p}}(\rd \bb{z})\bigg\|, \\
            \end{aligned}
        \end{equation}
        and the infima are taken over all Markov kernels $T_1 : \mathbb{K}_d \times \mathscr{B}(\R^d) \to [0,1]$ and $T_2 : \R^d \times \mathscr{B}(\mathbb{K}_d) \to [0,1]$.
    \end{theorem}

    Now, consider the following multivariate normal experiments with independent components
    \begin{alignat*}{6}
        &\overline{\mathscr{Q}}
        &&\leqdef &&~\{\overline{\QQ}_{n,\bb{p}}\}_{\bb{p}\in \Theta_R}, \quad &&\overline{\QQ}_{n,\bb{p}} ~\text{is the measure induced by } \mathrm{Normal}_d(n \bb{p}, n \mathrm{diag}(\bb{p})), \\
        &\mathscr{Q}^{\star}\hspace{-0.5mm}
        &&\leqdef &&~\{\QQ_{n,\bb{p}}^{\star}\}_{\bb{p}\in \Theta_R}, \quad &&\QQ_{n,\bb{p}}^{\star} ~\text{is the measure induced by } \mathrm{Normal}_d(\hspace{-0.75mm}\sqrt{n \bb{p}}, \mathrm{diag}(\bb{1}/4)),
    \end{alignat*}
    where $\bb{1} \leqdef (1,1,\dots,1)^{\top}$,
    then \citet[Section~7]{MR1922539} showed, using a variance stabilizing transformation, that
    \begin{equation}\label{eq:LeCam.distance.indep.normals}
        \Delta(\mathscr{Q},\overline{\mathscr{Q}}) \leq C_R \, \sqrt{\frac{d}{n}} \qquad \text{and} \qquad \Delta(\overline{\mathscr{Q}},\mathscr{Q}^{\star}) \leq C_R \, \frac{d}{\sqrt{n}},
    \end{equation}
    with proper adjustments to the definition of the deficiencies in \eqref{eq:def:deficiency.one.sided}.

    \begin{corollary}\label{cor:main.theorem.consequence}
        With the same notation as in Theorem~\ref{thm:bound.deficiency.distance}, we have, for $N\geq n^3 / d^{\hspace{0.2mm}2}$,
        \begin{equation}
            \Delta(\mathscr{P},\overline{\mathscr{Q}}) \leq C_R \, \frac{d}{\sqrt{n}} \qquad \text{and} \qquad \Delta(\mathscr{P},\mathscr{Q}^{\star}) \leq C_R \, \frac{d}{\sqrt{n}},
        \end{equation}
        where $C_R$ is a positive constant that depends only on $R$.
    \end{corollary}

\section{Proofs}\label{sec:proofs}

    \begin{proof}[\bf Proof of Theorem~\ref{thm:p.k.expansion}]
        Throughout the proof, the parameter $n\in \N$ satisfies $n\leq \gamma N$ and the asymptotic expressions are valid as $N\to \infty$.
        Let $\bb{p}\in N^{-1} \, \N_0^d \cap \mathrm{Int}\hspace{0.2mm}(\mathcal{S}_d)$ and $\bb{k}\in \mathbb{K}_d$.
        Using Stirling's formula,
        \begin{equation}
            \log m! = \frac{1}{2} \log(2\pi) + (m + \tfrac{1}{2}) \log m - m + \frac{1}{12 m} + \OO(m^{-3}), \quad m\to \infty,
        \end{equation}
        see, e.g., \citet[p.257]{MR0167642}, and taking logarithms in \eqref{eq:hypergeometric.pmf} and \eqref{eq:multinomial.pmf}, we obtain
        \begin{align}\label{eq:thm:p.k.expansion.eq.begin.1}
            \log\left(\frac{P_{N,n,\bb{p}}(\bb{k})}{Q_{n,\bb{p}}(\bb{k})}\right)
            &= \sum_{i=1}^{d+1} \log (N p_i)! - \sum_{i=1}^{d+1} \log (N p_i - k_i)! + \log (N - n)! - \log N! - \sum_{i=1}^{d+1} k_i \log p_i \notag \\
            &= \sum_{i=1}^{d+1} (N p_i + \tfrac{1}{2}) \log (N p_i) - \sum_{i=1}^{d+1} (N p_i - k_i + \tfrac{1}{2}) \log (N p_i - k_i) \notag \\[-2.5mm]
            &\qquad+ (N - n + \tfrac{1}{2}) \log (N - n) - (N + \tfrac{1}{2}) \log N - \sum_{i=1}^{d+1} k_i \log p_i \notag \\[-1.5mm]
            &\qquad+ \frac{1}{12 N} \left[\sum_{i=1}^{d+1} \frac{1}{p_i} \left\{1 - \bigg(1 - \frac{k_i}{N p_i}\bigg)^{-1}\right\} + \bigg(1 - \frac{n}{N}\bigg)^{-1} - 1\right] \notag \\
            &\qquad+ \OO\left(\frac{1}{N^3} \left[\sum_{i=1}^{d+1} \frac{1}{p_i^3} \left\{1 + \bigg(1 - \frac{k_i}{N p_i}\bigg)^{-3}\right\} + \bigg(1 - \frac{n}{N}\bigg)^{-3} + 1\right]\right) \notag \\
            &= - \sum_{i=1}^{d+1} N p_i \, \bigg(1 - \frac{k_i}{N p_i}\bigg) \log \bigg(1 - \frac{k_i}{N p_i}\bigg) - \frac{1}{2} \sum_{i=1}^{d+1} \log \bigg(1 - \frac{k_i}{N p_i}\bigg) \notag \\
            &\qquad+ N \left(1 - \frac{n}{N}\right) \log \left(1 - \frac{n}{N}\right) + \frac{1}{2} \log \left(1 - \frac{n}{N}\right) \notag \\
            &\qquad+ \frac{1}{12 N} \left[\sum_{i=1}^{d+1} \frac{1}{p_i} \left\{1 - \bigg(1 - \frac{k_i}{N p_i}\bigg)^{-1}\right\} + \bigg(1 - \frac{n}{N}\bigg)^{-1} - 1\right] \\
            &\qquad+ \OO_{\gamma}\left(\sum_{i=1}^{d+1} \frac{1}{(N p_i)^3}\right).
        \end{align}
        By applying the following Taylor expansions, valid for $|x| \leq \gamma < 1$,
        \begin{equation}\label{eq:thm:p.k.expansion.eq.Taylor}
            \begin{aligned}
                (1 - x) \log (1 - x) &= - x + \frac{x^2}{2} + \frac{x^3}{6} + \OO\big((1 - \gamma)^{-3} |x|^4\big), \\
                \log (1 - x) &= - x - \frac{x^2}{2} + \OO\big((1 - \gamma)^{-3} |x|^3\big), \\[2mm]
                (1 - x)^{-1} &= 1 + x + \OO\big((1 - \gamma)^{-3} |x|^2\big),
            \end{aligned}
        \end{equation}
        in \eqref{eq:thm:p.k.expansion.eq.begin.1}, we have
        \begin{equation}\label{eq:thm:p.k.expansion.eq.begin.2}
            \begin{aligned}
                \log\left(\frac{P_{N,n,\bb{p}}(\bb{k})}{Q_{n,\bb{p}}(\bb{k})}\right)
                &= \sum_{i=1}^{d+1} N p_i \,
                    \left\{\bigg(\frac{k_i}{N p_i}\bigg) - \frac{1}{2} \bigg(\frac{k_i}{N p_i}\bigg)^2 - \frac{1}{6} \bigg(\frac{k_i}{N p_i}\bigg)^3 + \OO_{\gamma}\left(\bigg|\frac{k_i}{N p_i}\bigg|^4\right)\right\} \\
                &\quad+ \frac{1}{2} \sum_{i=1}^{d+1}
                    \left\{\bigg(\frac{k_i}{N p_i}\bigg) + \frac{1}{2} \bigg(\frac{k_i}{N p_i}\bigg)^2 + \OO_{\gamma}\left(\bigg|\frac{k_i}{N p_i}\bigg|^3\right)\right\} \\
                &\quad+ N \left\{- \bigg(\frac{n}{N}\bigg) + \frac{1}{2} \bigg(\frac{n}{N}\bigg)^2 + \frac{1}{6} \bigg(\frac{n}{N}\bigg)^3 + \OO_{\gamma}\left(\bigg|\frac{n}{N}\bigg|^4\right)\right\} \\
                &\quad- \frac{1}{2} \left\{\bigg(\frac{n}{N}\bigg) + \frac{1}{2} \bigg(\frac{n}{N}\bigg)^2 + \OO_{\gamma}\left(\bigg|\frac{n}{N}\bigg|^3\right)\right\} \\
                &\quad+ \frac{1}{12 N} \left[\sum_{i=1}^{d+1} \frac{1}{p_i} \cdot \left\{- \frac{k_i}{N p_i} + \OO_{\gamma}\left(\bigg|\frac{k_i}{N p_i}\bigg|^2\right)\right\} + \left\{\frac{n}{N} + \OO_{\gamma}\left(\bigg|\frac{n}{N}\bigg|^2\right)\right\}\right] \\
                &\quad+ \OO_{\gamma}\left(\sum_{i=1}^{d+1} \frac{1}{(N p_i)^3}\right).
            \end{aligned}
        \end{equation}
        After rearranging some terms and noticing that $\sum_{i=1}^{d+1} k_i = n$, we get
        \begin{equation}\label{eq:thm:p.k.expansion.eq.end}
            \begin{aligned}
                \log\left(\frac{P_{N,n,\bb{p}}(\bb{k})}{Q_{n,\bb{p}}(\bb{k})}\right)
                &= \frac{1}{N} \left[\bigg(\frac{n^2}{2} - \frac{n}{2}\bigg) - \sum_{i=1}^{d+1} \frac{1}{p_i} \cdot \bigg(\frac{k_i^2}{2} - \frac{k_i}{2}\bigg)\right] \\
                &\quad+ \frac{1}{N^2} \left[\bigg(\frac{n^3}{6} - \frac{n^2}{4} + \frac{n}{12}\bigg) - \sum_{i=1}^{d+1} \frac{1}{p_i^2} \cdot \bigg(\frac{k_i^3}{6} - \frac{k_i^2}{4} + \frac{k_i}{12}\bigg)\right] \\
                &\quad+ \OO_{\gamma}\left(\frac{1}{N^3} \left[n^4 + \sum_{i=1}^{d+1} \frac{k_i^4}{p_i^3}\right]\right).
            \end{aligned}
        \end{equation}
        This proves \eqref{eq:thm:p.k.expansion.eq.2}.
        Equation~\eqref{eq:thm:p.k.expansion.eq.1} follows from the same arguments, simply by keeping fewer terms for the Taylor expansions in \eqref{eq:thm:p.k.expansion.eq.Taylor}. The details are omitted for conciseness.
    \end{proof}

    \begin{proof}[\bf Proof of Theorem~\ref{thm:total.variation.thm}]
        Define
        \begin{equation}
            A_{N,n,\bb{p}}(\gamma) \leqdef \left\{\bb{k}\in \mathbb{K}_d : \max_{1 \leq i \leq d + 1} \frac{k_i}{p_i} \leq \gamma N\right\} + \left(-\frac{1}{2}, \frac{1}{2}\right)^d, \quad \gamma > 0.
        \end{equation}
        By the comparison of the total variation norm with the Hellinger distance on page~726 of \citet{MR1922539}, we already know that
        \begin{equation}\label{eq:first.bound.total.variation}
            \|\widetilde{\PP}_{N,n,\bb{p}} - \widetilde{\QQ}_{n,\bb{p}}\| \leq \sqrt{2 \, \PP(\bb{X}\in A_{N,n,\bb{p}}^c(1/2)) + \EE\left[\log\bigg(\frac{\rd \widetilde{\PP}_{N,n,\bb{p}}}{\rd \widetilde{\QQ}_{n,\bb{p}}}(\bb{X})\bigg) \, \ind_{\{\bb{X}\in A_{N,n,\bb{p}}(1/2)\}}\right]}.
        \end{equation}
        By applying a union bound together with the large deviation bound for the (univariate) hypergeometric distribution in \citet[Equation~(4)]{MR3193733}, we get, for $N$ large enough,
        \begin{equation}\label{eq:concentration.bound}
            \begin{aligned}
                \PP(\bb{X}\in A_{N,n,\bb{p}}^c(1/2))
                &\leq \sum_{i=1}^{d+1} \PP\Bigl(K_i > \tfrac{N}{2 n} \cdot n p_i - 1\Bigr)
                \leq \sum_{i=1}^{d+1} \PP\Bigl(K_i > \nu_i \cdot n p_i\Bigr) \\
                &\leq \sum_{i=1}^{d+1} \left(\frac{1}{\nu_i}\right)^{n \nu_i p_i} \left(\frac{1 - p_i}{1 - \nu_i p_i}\right)^{n (1 - \nu_i p_i)},
            \end{aligned}
        \end{equation}
        where $\nu_i \leqdef \lceil p_i^{-1} - 1\rceil$, for all $1 \leq i \leq d + 1$.
        To estimate the expectation in \eqref{eq:first.bound.total.variation}, note that if $P_{N,n,\bb{p}}(\bb{x})$ and $Q_{n,\bb{p}}(\bb{x})$ denote the density functions associated with $\widetilde{\PP}_{N,n,\bb{p}}$ and $\widetilde{\QQ}_{n,\bb{p}}$ (i.e., $P_{N,n,\bb{p}}(\bb{x})$ is equal to $P_{N,n,\bb{p}}(\bb{k})$ whenever $\bb{k}\in \mathbb{K}_d$ is closest to $\bb{x}$, and analogously for $Q_{n,\bb{p}}(\bb{x})$), then, for $N$ large enough, we have
        \begin{equation}\label{eq:I.plus.II.plus.III}
            \begin{aligned}
                \left|\EE\left[\log\Bigg(\frac{\rd \widetilde{\PP}_{N,n,\bb{p}}}{\rd \widetilde{\QQ}_{n,\bb{p}}}(\bb{X})\Bigg) \, \ind_{\{\bb{X}\in A_{N,n,\bb{p}}(1/2)\}}\right]\right|
                &= \left|\EE\left[\log\bigg(\frac{P_{N,n,\bb{p}}(\bb{X})}{Q_{n,\bb{p}}(\bb{X})}\bigg) \, \ind_{\{\bb{X}\in A_{N,n,\bb{p}}(1/2)\}}\right]\right| \\
                &\leq \EE\left[\bigg|\log\bigg(\frac{P_{N,n,\bb{p}}(\bb{K})}{Q_{n,\bb{p}}(\bb{K})}\bigg)\bigg| \, \ind_{\{\bb{K}\in A_{N,n,\bb{p}}(3/4)\}}\right].
            \end{aligned}
        \end{equation}
        By Theorem~\ref{thm:p.k.expansion} with $\gamma = 3/4$, we find
        \begin{equation}
            \begin{aligned}
                \EE\left[\log\Bigg(\frac{\rd \widetilde{\PP}_{N,n,\bb{p}}}{\rd \widetilde{\QQ}_{n,\bb{p}}}(\bb{X})\Bigg) \, \ind_{\{\bb{X}\in A_{N,n,\bb{p}}(1/2)\}}\right]
                &= \OO\left(\frac{n^2}{N} + \sum_{i=1}^{d+1} \frac{\EE[K_i^2]}{N p_i}\right) \\
                &= \OO\left(\frac{n^2}{N} + \sum_{i=1}^{d+1} \frac{n^2 p_i^2}{N p_i}\right) = \OO\left(\frac{n^2}{N}\right).
            \end{aligned}
        \end{equation}
        Together with the large deviation bound in \eqref{eq:concentration.bound}, we deduce from \eqref{eq:first.bound.total.variation} that
        \begin{equation}\label{eq:total.variation.hypergeometric.multinomial}
            \|\widetilde{\PP}_{N,n,\bb{p}} - \widetilde{\QQ}_{n,\bb{p}}\| \leq \sqrt{2 \sum_{i=1}^{d+1} \left(\frac{1}{\nu_i}\right)^{n \nu_i p_i} \left(\frac{1 - p_i}{1 - \nu_i p_i}\right)^{n (1 - \nu_i p_i)} \hspace{-1mm}+ \OO\left(\frac{n^2}{N}\right)}.
        \end{equation}
        Also, by Lemma~3.1 in \cite{MR4249129} (a slightly weaker bound can be found in Lemma~2 of \cite{MR1922539}), we already know that
        \begin{equation}\label{eq:total.variation.multinomial.normal}
            \|\widetilde{\QQ}_{n,\bb{p}} - \QQ_{n,\bb{p}}\| = \OO\left(\frac{d}{\sqrt{n}} \sqrt{\frac{\max\{p_1,\dots,p_d,p_{d+1}\}}{\min\{p_1,\dots,p_d,p_{d+1}\}}}\right).
        \end{equation}
        Putting \eqref{eq:total.variation.hypergeometric.multinomial} and \eqref{eq:total.variation.multinomial.normal} together yields the conclusion.
    \end{proof}

    \begin{proof}[\bf Proof of Theorem~\ref{thm:bound.deficiency.distance}]
        By Theorem~\ref{thm:total.variation.thm} with our assumption $N\geq n^3 / d^{\hspace{0.2mm}2}$, we get the desired bound on $\delta(\mathscr{P},\mathscr{Q})$ by choosing the Markov kernel $T_1^{\star}$ that adds $\bb{U}\sim \mathrm{Uniform}\hspace{0.2mm}(-1/2,1/2)^d$ to $\bb{K}\sim \mathrm{Hypergeometric}\hspace{0.2mm}(N,n,\bb{p})$, namely
        \begin{equation}
            \begin{aligned}
                T_1^{\star}(\bb{k},B) \leqdef \int_{(-\frac{1}{2},\frac{1}{2})^d} \ind_{B}(\bb{k} + \bb{u}) \rd \bb{u}, \quad \bb{k}\in \mathbb{K}_d, ~B\in \mathscr{B}(\R^d).
            \end{aligned}
        \end{equation}
        To get the bound on $\delta(\mathscr{Q},\mathscr{P})$, it suffices to consider a Markov kernel $T_2^{\star}$ that inverts the effect of $T_1^{\star}$, i.e., rounding off every components of $\bb{Z}\sim \mathrm{Normal}_d(n \bb{p}, n \Sigma_{\bb{p}})$ to the nearest integer.
        Then, as explained by \citet[Section~5]{MR1922539}, we get
        \vspace{-1mm}
        \begin{align}
            \delta(\mathscr{Q},\mathscr{P})
            &\leq \bigg\|\PP_{N,n,\bb{p}} - \int_{\R^d} T_2^{\star}(\bb{z}, \cdot \, ) \, \QQ_{n,\bb{p}}(\rd \bb{z})\bigg\| \notag \\[0.5mm]
            &= \bigg\|\int_{\R^d} T_2^{\star}(\bb{z}, \cdot \, ) \int_{\mathbb{K}_d} T_1^{\star}(\bb{k}, \rd \bb{z}) \, \PP_{N,n,\bb{p}}(\rd \bb{k}) - \int_{\R^d} T_2^{\star}(\bb{z}, \cdot \, ) \, \QQ_{n,\bb{p}}(\rd \bb{z})\bigg\| \notag \\
            &\leq \bigg\|\int_{\mathbb{K}_d} T_1^{\star}(\bb{k}, \cdot \, ) \, \PP_{N,n,\bb{p}}(\rd \bb{k}) - \QQ_{n,\bb{p}}\bigg\|,
        \end{align}
        and we obtain the same bound by Theorem~\ref{thm:total.variation.thm}.
    \end{proof}

    \begin{proof}[\bf Proof of Corollary~\ref{cor:main.theorem.consequence}]
        This follows directly from Theorem~\ref{thm:bound.deficiency.distance}, Equation~\eqref{eq:LeCam.distance.indep.normals} and the fact that $\Delta(\cdot , \cdot)$ is a pseudometric (i.e., the triangle inequality is valid).
    \end{proof}

\section*{Acknowledgements}

\noindent
The author thanks the referee for his/her comments.

\vspace{3mm}
\noindent
{\bf Funding:} The author is supported by postdoctoral fellowships from the NSERC (PDF) and the FRQNT (B3X supplement and B3XR).

\section*{Declarations}

\noindent
{\bf Conflict of interest:} The author declares no conflict of interest.


%
%

\phantomsection
\addcontentsline{toc}{chapter}{References}

\bibliographystyle{authordate1}
\bibliography{Ouimet_2021_LLT_hypergeometric_bib}

\begin{thebibliography}{}

\bibitem[\protect\citename{Abramowitz \& Stegun, }1964]{MR0167642}
Abramowitz, M., \& Stegun, I.~A. 1964.
\newblock {\em Handbook of {M}athematical {F}unctions with {F}ormulas,
  {G}raphs, and {M}athematical {T}ables}.
\newblock National Bureau of Standards Applied Mathematics Series, vol. 55.
\newblock For sale by the Superintendent of Documents, U.S. Government Printing
  Office, Washington, D.C.
\newblock \href{http://www.ams.org/mathscinet-getitem?mr=MR0167642}{MR0167642}.

\bibitem[\protect\citename{Bhattacharya \& Ranga~Rao, }1976]{MR0436272}
Bhattacharya, R.~N., \& Ranga~Rao, R. 1976.
\newblock {\em Normal {A}pproximation and {A}symptotic {E}xpansions}.
\newblock John Wiley \& Sons, New York-London-Sydney.
\newblock \href{http://www.ams.org/mathscinet-getitem?mr=MR0436272}{MR0436272}.

\bibitem[\protect\citename{Brown {\em et~al.}, }2004]{MR2102503}
Brown, L.~D., Carter, A.~V., Low, M.~G., \& Zhang, C.-H. 2004.
\newblock Equivalence theory for density estimation, {P}oisson processes and
  {G}aussian white noise with drift.
\newblock {\em Ann. Statist.}, {\bf 32}(5), 2074--2097.
\newblock \href{http://www.ams.org/mathscinet-getitem?mr=MR2102503}{MR2102503}.

\bibitem[\protect\citename{Carter, }2002]{MR1922539}
Carter, A.~V. 2002.
\newblock Deficiency distance between multinomial and multivariate normal
  experiments. Dedicated to the memory of Lucien Le Cam.
\newblock {\em Ann. Statist.}, {\bf 30}(3), 708--730.
\newblock \href{http://www.ams.org/mathscinet-getitem?mr=MR1922539}{MR1922539}.

\bibitem[\protect\citename{Govindarajulu, }1965]{MR207011}
Govindarajulu, Z. 1965.
\newblock Normal approximations to the classical discrete distributions.
\newblock {\em Sankhy\={a} Ser. A}, {\bf 27}, 143--172.
\newblock \href{http://www.ams.org/mathscinet-getitem?mr=MR207011}{MR207011}.

\bibitem[\protect\citename{Johnson {\em et~al.}, }1997]{MR1429617}
Johnson, N.~L., Kotz, S., \& Balakrishnan, N. 1997.
\newblock {\em Discrete {M}ultivariate {D}istributions}.
\newblock Wiley Series in Probability and Statistics: Applied Probability and
  Statistics.
\newblock John Wiley \& Sons, Inc., New York.
\newblock \href{http://www.ams.org/mathscinet-getitem?mr=MR1429617}{MR1429617}.

\bibitem[\protect\citename{Kolassa, }1994]{MR1295242}
Kolassa, J.~E. 1994.
\newblock {\em Series {A}pproximation {M}ethods in {S}tatistics}.
\newblock Lecture Notes in Statistics, vol. 88.
\newblock Springer-Verlag, New York.
\newblock \href{http://www.ams.org/mathscinet-getitem?mr=MR1295242}{MR1295242}.

\bibitem[\protect\citename{Le~Cam \& Yang, }2000]{MR1784901}
Le~Cam, L., \& Yang, G.~L. 2000.
\newblock {\em Asymptotics in {S}tatistics}. Second edn.
\newblock Springer Series in Statistics.
\newblock Springer-Verlag, New York.
\newblock \href{http://www.ams.org/mathscinet-getitem?mr=MR1784901}{MR1784901}.

\bibitem[\protect\citename{Luh \& Pippenger, }2014]{MR3193733}
Luh, K., \& Pippenger, N. 2014.
\newblock Large-deviation bounds for sampling without replacement.
\newblock {\em Amer. Math. Monthly}, {\bf 121}(5), 449--454.
\newblock \href{http://www.ams.org/mathscinet-getitem?mr=MR3193733}{MR3193733}.

\bibitem[\protect\citename{Mariucci, }2016]{MR3850766}
Mariucci, E. 2016.
\newblock {L}e {C}am theory on the comparison of statistical models.
\newblock {\em Grad. J. Math.}, {\bf 1}(2), 81--91.
\newblock \href{http://www.ams.org/mathscinet-getitem?mr=MR3850766}{MR3850766}.

\bibitem[\protect\citename{Mattner \& Schulz, }2018]{MR3717995}
Mattner, L., \& Schulz, J. 2018.
\newblock On normal approximations to symmetric hypergeometric laws.
\newblock {\em Trans. Amer. Math. Soc.}, {\bf 370}(1), 727--748.
\newblock \href{http://www.ams.org/mathscinet-getitem?mr=MR3717995}{MR3717995}.

\bibitem[\protect\citename{Morgenstern, }1968]{MR0254884}
Morgenstern, D. 1968.
\newblock {\em Einf\"{u}hrung in die {W}ahrscheinlichkeitsrechnung und
  mathematische {S}tatistik [{I}n {G}erman]}.
\newblock Die Grundlehren der mathematischen Wissenschaften, Band 124.
\newblock Springer-Verlag, Berlin-New York.
\newblock \href{http://www.ams.org/mathscinet-getitem?mr=MR0254884}{MR0254884}.

\bibitem[\protect\citename{Nussbaum, }1996]{MR1425959}
Nussbaum, M. 1996.
\newblock Asymptotic equivalence of density estimation and {G}aussian white
  noise.
\newblock {\em Ann. Statist.}, {\bf 24}(6), 2399--2430.
\newblock \href{http://www.ams.org/mathscinet-getitem?mr=MR1425959}{MR1425959}.

\bibitem[\protect\citename{Ouimet, }2021]{MR4249129}
Ouimet, F. 2021.
\newblock A precise local limit theorem for the multinomial distribution and
  some applications.
\newblock {\em J. Statist. Plann. Inference}, {\bf 215}, 218--233.
\newblock \href{http://www.ams.org/mathscinet-getitem?mr=MR4249129}{MR4249129}.

\bibitem[\protect\citename{Prokhorov, }1953]{MR56861}
Prokhorov, Y.~V. 1953.
\newblock Asymptotic behavior of the binomial distribution.
\newblock {\em Uspekhi Mat. Nauk}, {\bf 8}(3(55)), 135--142.
\newblock \href{http://www.ams.org/mathscinet-getitem?mr=MR56861}{MR56861}.

\end{thebibliography}

\end{document}